\documentclass[12pt]{amsart}
\usepackage{tipa}
\usepackage{amssymb}
\usepackage{graphicx}
\theoremstyle{definition}

\theoremstyle{remark}

\numberwithin{equation}{section}

\begin{document}
\title[]{{\Large \bf Growth estimates for a class of subharmonic Functions
in a Half Space  $^{\ast}$ } }
\author{   Pan Guoshuang$^{1,2}$ and Deng Guantie$^{1, \ast\ast}$}%
\address{$^{1}$ Sch. Math. Sci. \& Lab. Math. Com. Sys. \\
   Beijing Normal University\\
    100875  Beijing, The People's Republic of China
    }%
\address{$^{2}$ Department of Public Basic Courses \\
   Beijing  Institute of Fashion and Technology \\
    100029  Beijing, The People's Republic of China
    }%
\email{denggt@bnu.edu.cn}%

\thanks{$\ast $ The project supported by NSFC (Grant No.10671022) and by RFDP (Grant
No.20060027023)
\endgraf $\ast\ast$:Corresponding author.}

\begin{abstract}
 \ \ {  A class of subharmonic  functions represented by
 the modified  kernels are proved to have the  growth estimates
 $u(x)= o(x_n^{1-\alpha}|x|^{m+\alpha})$ at infinity in the upper
 half space of ${\bf R}^{n}$, which generalizes the growth properties
 of analytic functions and harmonic functions. }
 \end{abstract}
\maketitle


 \section*{ 1. Introduction and Main Theorem}

\vspace{0.3cm}

  { Let ${\bf R}^{n} (n\geq3)$  denote the  $n$-dimensional Euclidean
space with points $x=(x_1,x_2,\cdots,x_{n-1},x_{n})=(x',x_n)$, where
$x' \in {\bf R}^{n-1}$ and $x_{n} \in {\bf R}$.  The boundary and
closure of an open  $\Omega$ of ${\bf R}^{n}$ are denoted by
$\partial{\Omega}$
 and $\overline{\Omega}$ respectively.
 The upper half-space $H$ is the set
 $H=\{x=(x',x_n)\in {\bf R}^{n}:\; x_n>0\}$, whose boundary is
 $\partial{H}$ .
    We  write $B(x,\rho)$ and $\partial B(x,\rho) $ for the open ball
    and the sphere of radius $\rho$  centered at $x$ in ${\bf R}^{n}$.
 We identify ${\bf R}^{n}$ with ${\bf R}^{n-1}\times {\bf R}$ and
${\bf R}^{n-1}$ with $ {\bf R}^{n-1}\times \{0\}$,
  with this convention we then have $ \partial {H}={\bf R}^{n-1}$,
  writing typical points $x,\ y \in {\bf R}^{n}$ as $x=(x',x_n),\
y=(y',y_n),$  where $x'=(x_1,x_2,\cdots,x_{n-1}),\
y'=(y_1,y_2,\cdots y_{n-1}) \in {\bf R}^{n-1}$ and putting
$$
x\cdot y=\sum_{j=1}^{n}x_jy_j=x'\cdot y'+x_ny_n,\ \ |x|=\sqrt{x\cdot
x},\ \ |x'|=\sqrt{x'\cdot x'}.
$$

  For $x\in{\bf R}^{n}\backslash\{0\}$, let([\textbf{16}])
$$
 E(x)=-r_n|x|^{2-n},
$$
where $|x|$ is the Euclidean norm, $r_n=\frac{1}{(n-2)\omega_{n}}$
and $\omega_{n}=\frac{2\pi^{\frac{n}{2}}}{\Gamma(\frac{n}{2})}$ is
the surface area of the unit sphere in ${\bf R}^{n} $. We know that
$E$ is locally integrable in ${\bf R}^{n} $.

 The Green function $G(x,y)$ for the upper half space
 $H$ is given by([\textbf{16}])
$$
 G(x,y)=E(x-y)-E(x-y^{\ast}) \qquad x,y\in\overline{H} ,\  x\neq y,
$$
where $^{\ast}$ denotes reflection in the boundary plane $\partial
H$ just as $y^{\ast}=(y_1,y_2,\cdots,y_{n-1},-y_n)$, then we define
the Poisson kernel $P(x,y')$ when $x\in H$ and $y'\in
\partial H $ by
$$
 P(x,y')=-\frac{\partial G(x,y)}{\partial
 y_n}\bigg|_{y_n=0}=\frac{2x_n}{\omega_n|x-(y',0)|^n}.
$$

  The  Dirichlet problem of upper half space is to find a function
 $u$ satisfying
$$
 u\in C^2(H), \eqno{(1.1)}
$$
$$
 \Delta u=0,   x\in H, \eqno{(1.2)}
$$
$$
 \lim_{x\rightarrow x'}u(x)=f(x')\ {\rm nontangentially  \  a.e.}x'\in \partial H, \eqno{(1.3)}
$$
where $f$ is a measurable function of ${\bf R}^{n-1} $. The Poisson
integral of the upper half space is defined by
$$
u(x)=P[f](x)=\int_{{\bf R}^{n-1}}P(x,y')f(y')dy'.
$$
 As we all know, the Poisson integral $P[f]$ exists if
$$
\int_{{\bf R}^{n-1}}\frac{|f(y')|}{1+|y'|^n} dy'<\infty.
$$
(see [\textbf{17,18}] and [\textbf{20}])In this paper, we will
consider measurable functions $f$ in ${\bf R}^{n-1}$ satisfying
$$
\int_{{\bf R}^{n-1}}\frac{|f(y')|}{1+|y'|^{n+m}}
dy'<\infty.\eqno{(1.4)}
$$
It is well known that the Poisson kernel $P(x,y')$ has a series
expansion in terms of the ultraspherical ( or Gegenbauer )
polynomials $ C^{\lambda}_{k}(t)\ (\lambda
=\frac{n}{2})$([\textbf{7}] and [\textbf{12}]). The latter can be
defined by a generating function
$$
(1-2tr+r^2)^{-\lambda } = \sum _{k=0}^{\infty}
C_{k}^{\lambda}(t)r^k, \eqno{(1.5)}
$$
where $|r|<1$, $ |t|\leq 1$ and $ \lambda > 0$. The coefficients $
C^{\lambda}_{k}(t) $ is called the ultraspherical ( or Gegenbauer )
polynomial of degree $ k $  associated with $ \lambda $, the
function $ C^{\lambda}_{k}(t) $ is a  polynomial of degree $ k $ in
$ t $.
 To obtain a solution of  Dirichlet problem for the boundary date $f$, as in
 [\textbf{3,5,15}] and [\textbf{20}], we use the following modified functions defined by
$$
 E_m(x-y)=\left\{\begin{array}{ll}
 E(x-y)  &   \mbox{when }   |y|\leq 1,  \\
 E(x-y) +\sum_{k=0}^{m-1}\frac{r_n|x|^k}{|y|^{n-2+k}}C^{\frac{n-2}{2}}_{k}
 \left( \frac{x\cdot y}{|x||y|}\right ) &
\mbox{when}\   |y|> 1.
 \end{array}\right.
$$
Then we can define modified Green function $G_m(x,y)$ and the
modified Poisson
 kernel $P_m(x,y')$ by([\textbf{1,2,4,11}] and [\textbf{20}])
$$
 G_m(x,y)=E_{m+1}(x-y)-E_{m+1}(x-y^{\ast}) \qquad x,y\in\overline{H}, \ x\neq
 y;\eqno{(1.6)}
$$
$$
 P_m(x,y')=\left\{\begin{array}{ll}
 P(x,y')  &   \mbox{when }   |y'|\leq 1  ,\\
 P(x,y') - \sum_{k=0}^{m-1}\frac{2x_n|x|^k}{\omega_n|y'|^{n+k}}C^{n/2}_{k}
 \left( \frac{x\cdot (y',0)}{|x||y'|}\right )&
\mbox{when}\   |y'|> 1.
 \end{array}\right.\eqno{(1.7)}
$$

   Siegel-Talvila([\textbf{3}]) have proved the
  following result:

\vspace{0.2cm}
 \noindent
{\bf Theorem A } Let $f$ be a measurable function in ${\bf R}^{n-1}$
satisfying (1.4), then the harmonic function
$$
v(x)= \int_{{\bf R}^{n-1}}P_m(x,y')f(y')dy' \quad x\in H
 \eqno{(1.8)}
$$
satisfies (1.1), (1.2), (1.3) and
$$
v(x)= o(x_n^{1-n}|x|^{m+n})  \quad  {\rm as}  \
|x|\rightarrow\infty.  \eqno{(1.9)}
$$
where $P_m(x,y')$ is defined by (1.7).

  In order to describe the asymptotic behaviour of subharmonic functions in
  half-spaces([\textbf{8,9}] and [\textbf{10}]),
  we establish the following theorems.

\vspace{0.2cm}
 \noindent
{\bf Theorem 1} Let $f$ be a measurable function in ${\bf R}^{n-1}$
satisfying (1.4), and $0< \alpha\leq n$. Let $v(x)$ be the harmonic
function defined by
 (1.8).
Then there exists $x_j\in H,\ \rho_j>0,$ such that
$$
\sum
_{j=1}^{\infty}\frac{\rho_j^{n-\alpha}}{|x_j|^{n-\alpha}}<\infty
\eqno{(1.10)}
$$
holds and
$$
v(x)= o(x_n^{1-\alpha}|x|^{m+\alpha})  \quad  {\rm as}  \
|x|\rightarrow\infty   \eqno{(1.11)}
 $$
holds in $H-G$. where $ G=\bigcup_{j=1}^\infty B(x_j,\rho_j)$.

\vspace{0.2cm}
 \noindent
 {\bf Remark } If $\alpha=n$, then (1.10) is a finite sum,
 the set $G$ is the union of finite balls, so (1.9) holds in $H$.
 This is just the result of Siegel-Talvila, therefore,
 our result (1.11) is the generalization of Theorem A.

  Next, we will generalize Theorem 1 to subharmonic functions.

\vspace{0.2cm}
 \noindent
{\bf Theorem 2 } Let $f$ be a measurable function in ${\bf R}^{n-1}$
satisfying (1.4), let $\mu$ be a positive  Borel  measure satisfying
$$
\int_H\frac{y_n}{1+|y|^{n+m}} d\mu(y)<\infty. \eqno{(1.12)}
$$
Write the subharmonic function
$$
u(x)= v(x)+h(x), \quad x\in H
$$
where $v(x)$ is the harmonic function defined by (1.8), $h(x)$ is
defined by
$$
h(x)= \int_H G_m(x,y)d\mu(y)
$$
and $G_m(x,y)$ is defined by (1.6). Then there exists $x_j\in H,\
\rho_j>0,$ such that (1.10) holds and
$$
u(x)= o(x_n^{1-\alpha}|x|^{m+\alpha})  \quad  {\rm as}  \
|x|\rightarrow\infty
 $$
holds in $H-G$, where $ G=\bigcup_{j=1}^\infty B(x_j,\rho_j)$ and
$0< \alpha<2$.

  Next we are concerned with minimal thinness at
infinity for $v(x)$ and $h(x)$, for a set $E \subset H$ and an open
set $F \subset {\bf R}^{n-1}$, we consider the capacity
$$
C(E;F)= \inf\int_{{\bf R}^{n-1}}g(y')dy'
$$
where the infimum is taken over all nonnegative measurable functions
$g$ such that $g=0$ outside $F$ and
$$
\int_{{\bf R}^{n-1}}\frac{g(y')}{|x-(y',0)|^n}dy'\geq 1 \qquad  {\rm
for \ all}\ x\in E.
$$
We say that $E \subset H$ is minimally thin at infinity if
$$
\sum _{i=1}^{\infty}2^{-in}C(E_i;F_i)<\infty,
$$
where $E_i=\{x\in E:\; 2^i\leq |x|<2^{i+1}\}$ and $F_i=\{x\in {\bf
R}^{n-1}:\; 2^i< |x|<2^{i+3}\}$.

\vspace{0.2cm}
 \noindent
{\bf Theorem 3 } Let $f$ be a measurable function in ${\bf R}^{n-1}$
satisfying (1.4), then there exists a set $E \subset H$ such that
$E$ is minimally thin at infinity and
$$
 \lim_{|x|\rightarrow\infty,x\in H-E}\frac{v(x)}{x_n|x|^m}=0.
$$

  Similarly, for $h(x)$, we can also conclude the following:

\vspace{0.2cm}
 \noindent
{\bf Corollary 1 } Let $\mu$ be a positive  Borel  measure
satisfying (1.12), then there exists a set $E \subset H$ such that
$E$ is minimally thin at infinity and
$$
 \lim_{|x|\rightarrow\infty,x\in H-E}\frac{h(x)}{x_n|x|^m}=0.
$$

Finally we are concerned with rarefiedness at infinity for $v(x)$
and $h(x)$, for a set $E \subset H$ and an open set $F \subset H$,
we consider the capacity
$$
C(E;F)= \inf\int_{H}g(y)d\mu(y)
$$
where the infimum is taken over all nonnegative measurable functions
$g$ such that $g=0$ outside $F$ and
$$
\int_{H}\frac{g(y)}{|x-y|^{n-1}}d\mu(y)\geq 1 \qquad  {\rm for \
all}\  x\in E.
$$
We say that $E \subset H$ is rarefied at infinity if
$$
\sum _{i=1}^{\infty}2^{-i(n-1)}C(E_i;F_i)<\infty,
$$
where $E_i$ is as in Theorem 3 and $F_i=\{x\in H:\; 2^i<
|x|<2^{i+3}\}$.

\vspace{0.2cm}
 \noindent
{\bf Theorem 4 } Let $\mu$ be a positive  Borel  measure satisfying
(1.12), then there exists a set $E \subset H$ such that $E$ is
rarefied at infinity and
$$
 \lim_{|x|\rightarrow\infty,x\in H-E}\frac{h(x)}{|x|^{m+1}}=0.
$$

  Similarly, for $v(x)$, we can also conclude the following:

\vspace{0.2cm}
 \noindent
{\bf Corollary 2 } Let $f$ be a measurable function in ${\bf
R}^{n-1}$ satisfying (1.4), then there exists a set $E \subset H$
such that $E$ is rarefied at infinity and
$$
 \lim_{|x|\rightarrow\infty,x\in H-E}\frac{v(x)}{|x|^{m+1}}=0.
$$

\vspace{0.4cm}

\section*{2.   Proof of Theorem }

\vspace{0.3cm}

Let $\mu$ be a positive Borel measure  in ${\bf R}^n,\ \beta\geq0$,
the maximal function $M(d\mu)(x)$ of order $\beta$ is defined by
$$
M(d\mu)(x)=\sup_{ 0<r<\infty}\frac{\mu(B(x,r))}{r^\beta},
$$
then the maximal function $M(d\mu)(x):{\bf R}^n \rightarrow
[0,\infty)$ is lower semicontinuous, hence measurable. To see this,
for any $ \lambda >0 $, let $D(\lambda)=\{x\in{\bf
R}^{n}:M(d\mu)(x)>\lambda\}$. Fix $x \in D(\lambda)$, then there
exists
 $r>0$ such that $\mu(B(x,r))>tr^\beta$ for some $t>\lambda$, and
there exists $ \delta>0$ satisfying
$(r+\delta)^\beta<\frac{tr^\beta}{\lambda}$. If $|y-x|<\delta$, then
$B(y,r+\delta)\supset B(x,r)$, therefore $\mu(B(y,r+\delta))\geq
tr^\beta >\lambda(r+\delta)^\beta$. Thus $B(x,\delta)\subset
D(\lambda)$. This proves that $D(\lambda)$ is open for each
$\lambda>0$.

 In order to obtain the results, we
need these lemmas below:

\vspace{0.2cm}
 \noindent
{\bf Lemma 1 } Let $\mu$ be a positive Borel measure  in ${\bf
R}^n,\ \beta\geq0,\ \mu({\bf R}^n)<\infty,$ for any $ \lambda \geq
5^{\beta} \mu({\bf R}^n)$, set
$$
E(\lambda)=\{x\in{\bf R}^{n}:|x|\geq2,M(d\mu)(x) >
\frac{\lambda}{|x|^{\beta}}\}
$$
then  there exists $ x_j\in E(\lambda)\  ,\ \rho_j> 0,\
j=1,2,\cdots$, such that
$$
E(\lambda) \subset \bigcup_{j=1}^\infty B(x_j,\rho_j) \eqno{(2.1)}
$$
and
$$
\sum _{j=1}^{\infty}\frac{\rho_j^{\beta}}{|x_j|^{\beta}}\leq
\frac{3\mu({\bf R}^n)5^{\beta}}{\lambda} .\eqno{(2.2)}
$$
Proof: Let $E_k(\lambda)=\{x\in E(\lambda):2^k\leq |x|<2^{k+1}\}$,
then  for any $ x \in E_k(\lambda),$ there exists $ r(x)>0$, such
that $\mu(B(x,r(x))) >\lambda(\frac{r(x)}{|x|})^{\beta} $, therefore
$r(x)\leq 2^{k-1}$.
 Since $E_k(\lambda)$ can be covered
by
 the union of a family of balls $\{B(x,r(x)):x\in E_k(\lambda) \}$,
 by the Vitali Lemma([\textbf{6}]), there exists $  \Lambda_k\subset E_k(\lambda)$,
$\Lambda_k$ is at most countable, such that $\{B(x,r(x)):x\in
\Lambda_k \}$ are disjoint and
$$
E_k(\lambda) \subset
 \cup_{x\in \Lambda_k} B(x,5r(x)),
$$
so
$$
E(\lambda)=\cup_{k=1}^\infty E_k(\lambda) \subset \cup_{k=1}^\infty
\cup_{x\in \Lambda_k} B(x,5r(x)).
$$

  On the other hand, note that $ \cup_{x\in \Lambda_k} B(x,r(x)) \subset \{x:2^{k-1}\leq
|x|<2^{k+2}\} $, so that
$$
 \sum_{x \in \Lambda_k}\frac{(5r(x))^{\beta}}{|x|^{\beta}}
\leq 5^\beta\sum_{x\in\Lambda_k}\frac{\mu(B(x,r(x)))}{\lambda} \leq
\frac{5^\beta}{\lambda} \mu\{x:2^{k-1}\leq |x|<2^{k+2}\}.
$$
Hence we obtain
$$
 \sum _{k=1}^{\infty}\sum
_{x \in \Lambda_k}\frac{(5r(x))^{\beta}}{|x|^{\beta}}
 \leq
 \sum _{k=1}^{\infty}\frac{5^\beta}{\lambda} \mu\{x:2^{k-1}\leq |x|<2^{k+2}\}
 \leq
\frac{3\mu({\bf R}^n)5^{\beta}}{\lambda}.
$$
  Rearrange $ \{x:x \in \Lambda_k,k=1,2,\cdots\} $ and $
\{5r(x):x \in \Lambda_k,k=1,2,\cdots\}
 $, we get $\{x_j\}$
and $\{\rho_j\}$ such that
 (2.1) and
(2.2) hold.

\vspace{0.2cm}
 \noindent
{\bf Lemma 2 }  Gegenbauer polynomials have the following
properties:\\
$(1)\  |C_k^{\lambda }(t)|\leq C_k^{\lambda }(1)=\frac{\Gamma
(2\lambda +k)}{\Gamma (2\lambda )\Gamma(k+1)}, \ \ |t|\leq 1 ;$\\
$(2)\  \frac{d}{dt}C_k^\lambda (t)=2\lambda C_{k-1}^{\lambda+1}(t),\ \ k \geq 1;$\\
$(3)\  \sum _{k=0}^{\infty} C_k^\lambda (1)r^k=(1-r)^{-2\lambda};$\\
$(4)\  |C^{\frac{n-2}{2}}_{k} \left( t\right )-C^{\frac{n-2}{2}}_{k}
\left( t^{\ast}\right )| \leq(n-2)C^{n/2}_{k-1} \left( 1\right
)|t-t^{\ast}|,\ \ |t|\leq 1, \ \ |t^{\ast}|\leq 1$.\\
Proof: (1) and (2) can be derived from [\textbf{7}] and
[\textbf{13}]; (3) follows by taking $t=1$ in (1.5); (4) follows
 by (1), (2) and the Mean Value Theorem for Derivatives.

\vspace{0.2cm}
 \noindent
{\bf Lemma 3 }  Green function $G(x,y)$ has the following
estimates:\\
$(1)\  |G(x,y)|\leq \frac{r_n}{|x-y|^{n-2}};$\\
$(2)\  |G(x,y)|\leq \frac{2x_ny_n}{\omega_n|x-y|^n};$\\
$(3)\  |G(x,y)|\leq \frac{Ax_ny_n}{|x-y|^{n-2}|x-y^{\ast}|^2}.$\\
Proof: (1) is obvious; (2) follows by the Mean Value Theorem for
Derivatives; (3) can be derived from [\textbf{14}].

  Throughout the paper, let $A$
denote various positive constants independent of the variables in
question.

 \emph{Proof of Theorem 1}

Define the measure $dm(y')$ and the kernel $K(x,y')$ by
$$
dm(y')=\frac{|f(y')|}{1+|y'|^{n+m}} dy' ,\ \ K(x,y')=
P_m(x,y')(1+|y'|^{n+m}).
$$
  For any $\varepsilon >0$, there exists $R_\varepsilon >2$, such that
$$
\int_{|y'|\geq
R_\varepsilon}dm(y')\leq\frac{\varepsilon}{5^{n-\alpha}}.
$$
For every Lebesgue measurable set $E \subset {\bf R}^{n-1}$  , the
measure $m^{(\varepsilon)}$ defined by $m^{(\varepsilon)}(E)
=m(E\cap\{x'\in{\bf R}^{n-1}:|x'|\geq R_\varepsilon\}) $ satisfies
$m^{(\varepsilon)}({\bf
R}^{n-1})\leq\frac{\varepsilon}{5^{n-\alpha}}$, write
\begin{eqnarray*}
&v_1(x)& =\int_{|x-(y',0)| \leq 3|x|} P(x,y')(1+|y'|^{n+m})
dm^{(\varepsilon)}(y'),\\
&v_2(x)&=\int_{|x-(y',0)| \leq 3|x|}
(P_m(x,y')-P(x,y'))(1+|y'|^{n+m})
dm^{(\varepsilon)}(y'), \\
&v_3(x)&=\int_{|x-(y',0)| > 3|x|} K(x,y')dm^{(\varepsilon)}(y'), \\
&v_4(x)&=\int_{1<|y'|<R_\varepsilon}K(x,y') dm(y'), \\
&v_5(x)&=\int_{|y'|\leq1}K(x,y') dm(y'). \\
\end{eqnarray*}
then
$$
|v(x)| \leq |v_1(x)|+|v_2(x)|+|v_3(x)|+|v_4(x)|+|v_5(x)|.
\eqno{(2.3)}
$$
Let $ E_1(\lambda)=\{x\in{\bf R}^{n}:|x|\geq2,\exists
t>0,m^{(\varepsilon)}(B(x,t)\cap{\bf R}^{n-1}
)>\lambda(\frac{t}{|x|})^{n-\alpha}\}$, therefore, if $ |x|\geq
2R_\varepsilon$ and $ x \notin E_1(\lambda)$, then we have

\begin{eqnarray*}
|v_1(x)|
&\leq& \int_{x_n\leq|x-(y',0)| \leq
3|x|}\frac{2x_n}{\omega_n|x-(y',0)|^n}2|y'|^{n+m} dm^{(\varepsilon)}(y') \\
&\leq& \frac{4^{n+m+1}}{\omega_n}x_n|x|^{m+n}\int_{x_n}^{3|x|}
\frac{1}{t^n} dm_x^{(\varepsilon)}(t) \\
&\leq& \frac{4^{n+m+1}}{\omega_n}
\bigg(\frac{1}{3^\alpha}+\frac{n}{\alpha}\bigg)\lambda
x_n^{1-\alpha}|x|^{m+\alpha}.\hspace{39mm} (2.4)
\end{eqnarray*}
where  $m_x^{(\varepsilon)}(t)=\int_{|x-(y',0)| \leq t}
dm^{(\varepsilon)}(y')$.

   By (1) and (3) of Lemma 2, we obtain
\begin{eqnarray*}
|v_2(x)|
&\leq& \int_{x_n\leq|x-(y',0)| \leq 3|x|}
\sum_{k=0}^{m-1}\frac{2x_n|x|^k}{\omega_n}C^{n/2}_{k}
 \left( 1\right )\frac{2|y'|^{n+m}}{|y'|^{n+k}} dm^{(\varepsilon)}(y') \\
&\leq&
\frac{4^{m+1}}{\omega_n}\sum_{k=0}^{m-1}\frac{1}{4^k}C^{n/2}_{k}
 \left( 1\right )\frac{1}{5^{n-\alpha}}\varepsilon x_n|x|^m \\
&\leq& \frac{4^{m+1+\alpha}}{\omega_n\cdot 3^n}\varepsilon x_n|x|^m.
\hspace{67mm} (2.5)
\end{eqnarray*}

  By (1) and (3) of Lemma 2, we see that([\textbf{19}])
\begin{eqnarray*}
|v_3(x)|
&\leq& \int_{|x-(y',0)| > 3|x|}
\sum_{k=m}^{\infty}\frac{4x_n|x|^k}{\omega_n(2|x|)^{k-m}}C^{n/2}_{k}
\left( 1\right ) dm^{(\varepsilon)}(y') \\
&\leq& \frac{2^{m+2}}{\omega_n}\frac{\varepsilon}{5^{n-\alpha}}
\sum_{k=m}^{\infty}\frac{1}{2^k}C^{n/2}_{k}
\left( 1\right )x_n|x|^m \\
&\leq& \frac{2^{m-n+2\alpha+2}}{\omega_n}\varepsilon x_n|x|^m .
\hspace{62mm} (2.6)
\end{eqnarray*}

  Write
\begin{eqnarray*}
v_4(x)
&=& \int_{1<|y'|<R_\varepsilon}[P(x,y') +(P_m(x,y')-P(x,y'))](1+|y'|^{n+m}) dm(y') \\
&=& v_{41}(x)+v_{42}(x),
\end{eqnarray*}
then
\begin{eqnarray*}
|v_{41}(x)|
&\leq& \int_{1<|y'|<R_\varepsilon}\frac{2x_n}{\omega_n|x-(y',0)|^n}
2|y'|^{n+m} dm(y') \\
&\leq& \frac{4R_\varepsilon^{n+m}x_n}{\omega_n}
\int_{1<|y'|<R_\varepsilon}\frac{1}{(\frac{|x|}{2})^n}
 dm(y') \\
&\leq& \frac{2^{n+2}R_\varepsilon^{n+m} m({\bf R}^{n-1})}{\omega_n}
\frac{x_n}{|x|^n}.\hspace{50mm} (2.7)
\end{eqnarray*}
by (1) and (3) of Lemma 2, we obtain
\begin{eqnarray*}
|v_{42}(x)|
&\leq& \int_{1<|y'|<R_\varepsilon}
\sum_{k=0}^{m-1}\frac{2x_n|x|^k}{\omega_n|y'|^{n+k}}C^{n/2}_{k}
 \left( 1\right )\cdot 2|y'|^{n+m} dm(y') \\
&\leq& \sum_{k=0}^{m-1}\frac{4}{\omega_n}C^{n/2}_{k} \left( 1\right
)x_n|x|^kR_\varepsilon^{m-k} m({\bf
R}^{n-1}) \\
&\leq& \frac{2^{n+m+1}R_\varepsilon^m m({\bf R}^{n-1})}{\omega_n}
x_n|x|^{m-1}.\hspace{41mm} (2.8)
\end{eqnarray*}
  In case $|y'|\leq 1$, note that
$$
K(x,y')=P_m(x,y')(1+|y'|^{n+m})
\leq\frac{4x_n}{\omega_n|x-(y',0)|^n},
$$
so that
$$
|v_5(x)|\leq \int_{|y'|\leq1}\frac{4x_n}{\omega_n(\frac{|x|}{2})^n}
dm(y') \leq\frac{2^{n+2} m({\bf R}^{n-1})}{\omega_n}
\frac{x_n}{|x|^n}.\eqno{(2.9)}
$$

  Thus, by collecting (2.3), (2.4), (2.5), (2.6), (2.7), (2.8) and
(2.9), there exists a positive constant $A$ independent of
$\varepsilon$, such that if $ |x|\geq 2R_\varepsilon$ and $\  x
\notin E_1(\varepsilon)$, we have
$$
|v(x)|\leq A\varepsilon x_n^{1-\alpha}|x|^{m+\alpha}.
$$

 Let $\mu_\varepsilon$ be a measure in ${\bf R}^n$ defined by
$ \mu_\varepsilon(E)= m^{(\varepsilon)}(E\cap{\bf R}^{n-1})$ for
every measurable set $E$ in ${\bf R}^n$.Take
$\varepsilon=\varepsilon_p=\frac{1}{2^{p+2}}, p=1,2,3,\cdots$, then
there exists a sequence $ \{R_p\}$: $1=R_0<R_1<R_2<\cdots$ such that
$$
\mu_{\varepsilon_p}({\bf R}^n)=\int_{|y'|\geq
R_p}dm(y')<\frac{\varepsilon_p}{5^{n-\alpha}}.
$$
Take $\lambda=3\cdot5^{n-\alpha}\cdot2^p\mu_{\varepsilon_p}({\bf
R}^n)$ in Lemma 1, then there exists $ x_{j,p}$ and $ \rho_{j,p}$,
where $R_{p-1}\leq |x_{j,p}|<R_p,$ such that
$$
\sum _{j=1}^{\infty}(\frac{\rho_{j,p}}{|x_{j,p}|})^{n-\alpha} \leq
\frac{1}{2^{p}}.
$$
if $R_{p-1}\leq |x|<R_p$ and $ x\notin G_p=\cup_{j=1}^\infty
B(x_{j,p},\rho_{j,p})$, we have
$$
|v(x)|\leq A\varepsilon_px_n^{1-\alpha}|x|^{m+\alpha},
$$
Thereby
$$
\sum _{p=1}^{\infty}
\sum_{j=1}^{\infty}(\frac{\rho_{j,p}}{|x_{j,p}|})^{n-\alpha} \leq
\sum _{p=1}^{\infty}\frac{1}{2^{p}}=1<\infty.
$$

  Set $ G=\cup_{p=1}^\infty G_p$, thus Theorem 1 holds.

 \emph{Proof of Theorem 2}

  Define the measure $dn(y)$ and the kernel $L(x,y)$ by
$$
dn(y)=\frac{y_n d\mu(y)}{1+|y|^{n+m}},\ \
L(x,y)=G_m(x,y)\frac{1+|y|^{n+m}}{y_n}.
$$
then the function $h(x)$ can be written as
$$
h(x)=\int_H L(x,y) dn(y).
$$

  For any $\varepsilon >0$, there exists $R_\varepsilon >2$, such that
$$
\int_{|y|\geq R_\varepsilon}dn(y)<\frac{\varepsilon}{5^{n-\alpha}}.
$$
For every Lebesgue measurable set $E \subset {\bf R}^{n}$, the
measure $n^{(\varepsilon)}$ defined by $n^{(\varepsilon)}(E)
=n(E\cap\{y\in H:|y|\geq R_\varepsilon\}) $ satisfies
$n^{(\varepsilon)}(H)\leq\frac{\varepsilon}{5^{n-\alpha}}$, write
\begin{eqnarray*}
&h_1(x)& =\int_{|x-y|\leq\frac{x_n}{2}}G(x,y)\frac{1+|y|^{n+m}}{y_n}
dn^{(\varepsilon)}(y), \\
&h_2(x)&=\int_{\frac{x_n}{2}<|x-y|\leq3|x|}G(x,y)\frac{1+|y|^{n+m}}{y_n}
dn^{(\varepsilon)}(y), \\
&h_3(x)&=\int_{|x-y|\leq3|x|}(G_m(x,y)-G(x,y))\frac{1+|y|^{n+m}}{y_n}
dn^{(\varepsilon)}(y), \\
&h_4(x)&=\int_{|x-y|>3|x|}L(x,y)
dn^{(\varepsilon)}(y), \\
&h_5(x)&=\int_{1<|y|<R_\varepsilon}L(x,y) dn(y), \\
&h_6(x)&=\int_{|y|\leq1}L(x,y) dn(y). \\
\end{eqnarray*}
then
$$
h(x)=h_1(x)+h_2(x)+h_3(x)+h_4(x)+h_5(x)+h_6(x). \eqno{(2.10)}
$$
Let $ E_2(\lambda)=\{x\in{\bf R}^{n}:|x|\geq2,\exists
t>0,n^{(\varepsilon)}(B(x,t)\cap H
)>\lambda(\frac{t}{|x|})^{n-\alpha}\}, $ therefore, if $ |x|\geq
2R_\varepsilon$ and $ x \notin E_1(\lambda)$, then we have by (1) of
Lemma 3
\begin{eqnarray*}
|h_1(x)|
&\leq& \int_{|x-y|\leq\frac{x_n}{2}}
\frac{r_n}{|x-y|^{n-2}}\frac{2|y|^{n+m}}{\frac{x_n}{2}}
dn^{(\varepsilon)}(y) \\
&\leq& 4\times
(3/2)^{n+m}r_n\frac{|x|^{n+m}}{x_n}\int_0^\frac{x_n}{2}
\frac{1}{t^{n-2}} dn_x^{(\varepsilon)}(t)\\
&\leq& 4\times (3/2)^{n+m}r_n\bigg[\frac{1}{2^{2-\alpha}}+
\frac{n-2}{(2-\alpha)2^{2-\alpha}}\bigg]\lambda
x_n^{1-\alpha}|x|^{m+\alpha}.\hspace{2mm} (2.11)
\end{eqnarray*}
where $ n_x^{(\varepsilon)}(t)=\int_{|x-y| \leq t}
dn^{(\varepsilon)}(y)$.\\
By (2) of Lemma 3, we have
\begin{eqnarray*}
|h_2(x)|
&\leq&
\int_{\frac{x_n}{2}<|x-y|\leq3|x|}\frac{2x_ny_n}{\omega_n|x-y|^n}
\frac{2|y|^{n+m}}{y_n}
dn^{(\varepsilon)}(y) \\
&\leq&
\frac{4^{n+m+1}}{\omega_n}x_n|x|^{n+m}\int_\frac{x_n}{2}^{3|x|}
\frac{1}{t^n} dn_x^{(\varepsilon)}(t)\\
&\leq& \frac{4^{n+m+1}}{\omega_n}\bigg(\frac{1}{3^\alpha}+
\frac{n2^\alpha}{\alpha}\bigg)\lambda x_n^{1-\alpha}|x|^{m+\alpha}.
\hspace{33mm} (2.12)
\end{eqnarray*}

  First note $C^{\lambda}_{0}
\left( t\right )\equiv 1([\textbf{7}])$ , then we obtain by (1), (3)
and (4) of Lemma 2 and taking $t=\frac{x\cdot y}{|x||y|},\
t^{\ast}=\frac{x\cdot y^{\ast}}{|x||y^{\ast}|}$ in (4) of Lemma 2

\begin{eqnarray*}
|h_3(x)|
&\leq&
\int_{|x-y|\leq3|x|}\sum_{k=1}^{m}\frac{r_n|x|^k}{|y|^{n-2+k}}
2(n-2)C^{n/2}_{k-1} \left( 1\right )\frac{x_ny_n}{|x||y|}
\frac{2|y|^{n+m}}{y_n} dn^{(\varepsilon)}(y) \\
&\leq&
\frac{4^{m+1}}{\omega_n}\sum_{k=1}^{m}\frac{1}{4^{k-1}}C^{n/2}_{k-1}
 \left( 1\right )\frac{1}{5^{n-\alpha}}\varepsilon x_n|x|^m \\
&\leq& \frac{4^{m+1+\alpha}}{\omega_n\cdot 3^n}\varepsilon x_n|x|^m.
\hspace{66mm} (2.13)
\end{eqnarray*}

  By (1), (3) and (4) of Lemma 2, we see that
\begin{eqnarray*}
|h_4(x)|
&\leq&
\int_{|x-y|>3|x|}\sum_{k=m+1}^{\infty}\frac{r_n|x|^k}{|y|^{n-2+k}}
2(n-2)C^{n/2}_{k-1} \left( 1\right )\frac{x_ny_n}{|x||y|}
\frac{2|y|^{n+m}}{y_n} dn^{(\varepsilon)}(y) \\
&\leq&
\frac{2^{m+2}}{\omega_n}\sum_{k=m+1}^{\infty}\frac{1}{2^{k-1}}C^{n/2}_{k-1}
 \left( 1\right )\frac{1}{5^{n-\alpha}}\varepsilon x_n|x|^m \\
&\leq& \frac{2^{m-n+2\alpha+2}}{\omega_n}\varepsilon x_n|x|^m .
\hspace{60mm} (2.14)
\end{eqnarray*}

  Write
\begin{eqnarray*}
h_5(x)
&=&
\int_{1<|y|<R_\varepsilon}[G(x,y)+(G_m(x,y)-G(x,y))]\frac{1+|y|^{n+m}}{y_n} dn(y) \\
&=& h_{51}(x)+h_{52}(x),
\end{eqnarray*}
then we obtain by (2) of Lemma 3
\begin{eqnarray*}
|h_{51}(x)|
&\leq& \int_{1<|y|<R_\varepsilon}
\frac{2x_ny_n}{\omega_n|x-y|^n}\frac{2|y|^{n+m}}{y_n} dn(y) \\
&\leq&
\frac{4R_\varepsilon^{n+m}}{\omega_n}x_n\int_{1<|y|<R_\varepsilon}
\frac{1}{(\frac{|x|}{2})^n} dn(y) \\
&\leq& \frac{2^{n+2}R_\varepsilon^{n+m}n(H) }{\omega_n}
\frac{x_n}{|x|^n} .\hspace{56mm} (2.15)
\end{eqnarray*}
by (1), (3) and (4) of Lemma 2, we obtain
\begin{eqnarray*}
|h_{52}(x)|
&\leq&
\int_{1<|y|<R_\varepsilon}\sum_{k=1}^{m}\frac{r_n|x|^k}{|y|^{n-2+k}}
2(n-2)C^{n/2}_{k-1} \left( 1\right )\frac{x_ny_n}{|x||y|}
\frac{2|y|^{n+m}}{y_n} dn(y) \\
&\leq& \sum_{k=1}^{m}\frac{4}{\omega_n}C^{n/2}_{k-1}
 \left( 1\right )x_n|x|^{k-1}R_\varepsilon^{m-k+1}n(H) \\
&\leq& \frac{2^{n+m+1}R_\varepsilon^{m}n(H) }{\omega_n} x_n|x|^{m-1}
.\hspace{47mm} (2.16)
\end{eqnarray*}

  In case $|y|\leq 1$, by (2) of Lemma 3, we have
$$
|L(x,y)|\leq \frac{2x_ny_n}{\omega_n|x-y|^n}\frac{2}{y_n}
=\frac{4x_n}{\omega_n|x-y|^n},
$$
so that
$$
|h_6(x)|\leq \int_{|y|\leq1}\frac{4x_n}{\omega_n(\frac{|x|}{2})^n} 
dn(y)\leq \frac{2^{n+2}n(H) }{\omega_n}
\frac{x_n}{|x|^n}.\eqno{(2.17)}
$$

  Thus, by collecting (2.10), (2.11), (2.12), (2.13),
(2.14), (2.15), (2.16) and (2.17), there exists a positive constant
$A$ independent of $\varepsilon$, such that if $ |x|\geq
2R_\varepsilon$ and $\  x \notin E_2(\varepsilon)$, we have
$$
 |h(x)|\leq A\varepsilon x_n^{1-\alpha}|x|^{m+\alpha}.
$$

  Similarly, if $x\notin G$, we have
$$
h(x)= o(x_n^{1-\alpha}|x|^{m+\alpha})\quad  {\rm as} \
|x|\rightarrow\infty. \eqno{(2.18)}
$$

by (1.11) and  (2.18), we obtain
$$
u(x)=v(x)+h(x)= o(x_n^{1-\alpha}|x|^{m+\alpha})\quad  {\rm as} \
|x|\rightarrow\infty
$$
hold in $H-G$, thus we complete the proof of Theorem 2.

 \emph{Proof of Theorem 3 and 4}

We prove only Theorem 4, the proof of Theorem 3 is similar. By
(2.13), (2.14), (2.15), (2.16) and (2.17) we have
$$
 \lim_{|x|\rightarrow\infty,x\in H}\frac{h_3(x)+h_4(x)+h_5(x)+h_6(x)}{|x|^{m+1}}=0.\eqno{(2.19)}
$$
In view of (1.12), we can find a sequence $\{a_i\}$ of positive
numbers such that $\lim_{i\rightarrow\infty}a_i=\infty$ and
$$
\sum
_{i=1}^{\infty}a_i\int_{F_i}\frac{y_n}{|y|^{n+m}}d\mu(y)<\infty.
$$
Consider the sets
$$
E_i=\{x\in H:\; 2^i\leq |x|<2^{i+1},|h_1(x)+h_2(x)|\geq
a_i^{-1}2^{im}|x|\}
$$
for $i=1,2,\cdots$. If $x\in E_i$, then we obtain by (3) of Lemma 3
$$
a_i^{-1}\leq 2^{-im}|x|^{-1}|h_1(x)+h_2(x)|\leq
A2^{-i(m+1)}\int_{F_i}\frac{y_n}{|x-y|^{n-1}}d\mu(y)
$$
so that it follows from the definition of $C(E_i;F_i)$ that
$$
C(E_i;F_i)\leq Aa_i2^{-i(m+1)}\int_{F_i}y_nd\mu(y)\leq
Aa_i2^{i(n-1)}\int_{F_i}\frac{y_n}{|y|^{n+m}}d\mu(y)
$$
Define $E=\bigcup_{i=1}^\infty E_i$, then
$$
\sum _{i=1}^{\infty}2^{-i(n-1)}C(E_i;F_i)<\infty.
$$
Clearly,
$$
 \lim_{|x|\rightarrow\infty,x\in H-E}\frac{h_1(x)+h_2(x)}{|x|^{m+1}}=0.\eqno{(2.20)}
$$
Thus, by collecting (2.19) and (2.20), the proof of Theorem 4 is
completed.

\begin{center}

\end{center}

\end{document}